\documentclass[graybox]{svmult}

\usepackage{type1cm}        
\usepackage{makeidx}         
\usepackage{graphicx}        
\usepackage{multicol}        
\usepackage[bottom]{footmisc}

\usepackage{newtxtext}       %
\usepackage[varvw]{newtxmath}       

\usepackage{latexsym}
\usepackage{graphicx}

\usepackage{color}
\usepackage{bm}

\def\ub {{\bm u}}

\def\0b {{\bf 0}}

\def\ub {{\bf u}}

\def\Wb {{\bf W}}
\def\Vb {{\bf V}}

\def\ub {\bm{u}}

\def\nb {\bm{n}}

\def\ub {\bm{u}}

\def\nb {\bm{n}}

\makeindex             

\begin{document}

\title*{FETI-DP algorithms for 2D Biot model with discontinuous Galerkin discretization}
\titlerunning{FETI-DP algorithm for 2D Biot model}
\author{Pilhwa Lee}
\institute{Department of Mathematics, Morgan State University, 1700 E. Cold Spring Lane, Baltimore, MD, USA
\email{Pilhwa.Lee@morgan.edu}}
%
%
\maketitle

\abstract{Dual-Primal Finite Element Tearing and Interconnecting (FETI-DP) algorithms are developed for a 2D Biot model. The model is formulated with mixed-finite elements as a saddle-point problem. The displacement $\mathbf{u}$ and the Darcy flux flow $\mathbf{z}$ are represented with $P_1$ piecewise continuous elements and pore-pressure $p$ with $P_0$ piecewise constant elements, {\it i.e.}, overall three fields with a stabilizing term. We have tested the functionality of FETI-DP with Dirichlet preconditioners. Numerical experiments show a signature of scalability of the resulting parallel algorithm in the compressible elasticity with permeable Darcy flow as well as almost incompressible elasticity.}

\section{Introduction}

\indent \indent Poroelasticity, {\it i.e.}, elasticity of porous media with permeated Darcy flow, pioneered by Biot \cite{Biot1941general, Biot1955theory} has been used broadly in geoscience \cite{Phillips2008} and biomechanics \cite{Badia2009, Chapelle2010, Richardson2021} among many others. The difficulties for solving the linear elasticity and incompressible flow problems also arise in solving the poroelastic problem, and there have been diverse mathematical formulations and discretizations. When a continuous Galerkin approach was formulated with mixed finite elements with three-fields of displacement, Darcy flow flux, and pressure \cite{Phillips2007}, the main numerical difficulties are elastic locking and non-physical oscillatory pressure profiles. There have been some new methods for dealing with these difficulties; for example, continuous Galerkin with non-standard three-fields of displacement, pressure, and volumetric stress \cite{Oyarzua2016locking}, discontinuous Galerkin formulations \cite{Phillips2008, Hong2018} with standard three-fields as well as non-conforming mixed finite elements \cite{Yi2011, Boffi2016}. When lowest-order finite elements are applied, stabilizing terms should be added \cite{Burman2007, Berger2015, Rodrigo2018} to satisfy the inf-sup condition \cite{Brezzi1991, Howell2011}.\\
\indent In this paper, we propose a numerical scheme for solving the Biot model with three-fields linear poroelasticity. We consider a discontinuous Galerkin discretization, {\it i.e.}, the displacement and Darcy flow flux discretized as piecewise continuous in $P_1$ elements, and the pore pressure as piecewise constant in the $P_0$ space with a stabilizing term. The emerging formulation is a saddle-point problem, and more specifically, a twofold saddle-point problem. This indefinite system is computationally challenging with slow convergence in iterative methods. It is necessary to incorporate relevant preconditioners for saddle-point problems \cite{Benzi2005numerical, Mardal2011preconditioning}.\\
\indent FETI-DP algorithms transform indefinite problems to positive definite interface problems of Lagrangian multipliers for subdomains and a primal problem for the coarse space \cite{Farhat2001}. They have been applied for linear elasticity \cite{Klawonn2006, Pavarino2016} and incompressible Stokes flows \cite{Li2005, Kim2006, Tu2013, Tu2015} as saddle-point problems. There are theoretical bounds for the condition numbers in the preconditioned systems independent to partitioned subdomains. We show numerical scalability of FETI-DP algorithm preconditioned by Dirichlet preconditioner for the three-fields Biot model discretized with stabilized $P_1-P_1-P_0$.

\section{Linear poroelastic model}
Poroelastic models describe the interaction of fluid flows and deformable elastic porous media saturated in the
fluid. Let $\bm{u}$ be the elastic displacement, 
$p$ be the pore-pressure. We assume that the permeability is homogeneous: ${\bf K} = \kappa{\bf I}$. Denote $\bm{z}$ as the Darcy volumetric fluid flux. The quasi-static Biot model reads as:
\begin{eqnarray}
-(\lambda + \mu) \nabla (\nabla \cdot \mathbf{u}) - \mu \nabla^2 \mathbf{u} + \alpha \nabla p &=& \mathbf{f},\\
{\bf K}^{-1} \bm{z} + \nabla p &=& \mathbf{b},   \label{lee2}\\
\frac{\partial}{\partial t} \left(\alpha \nabla \cdot \bm{u} + c_0 p \right)+ \nabla \cdot \bm{z} &=& g. \label{lee3}
\end{eqnarray}
The first equation is the moment conservation. The second equation is Darcy's law. The third equations is the mass conservation equation. For simplicity, we neglect the effects of gravity acceleration. In the above equations, $\bm{f}$ is the body force of the solid, $\bm{b}$ is the body force of the fluid, $g$ is a source or sink term, $c_0>0$ is the constrained specific storage coefficient, $\alpha$ is the Biot-Willis constant which is close to 1. $\lambda$ and $\mu$ are the first and second Lam\'e parameters, respectively.


We consider $\Omega \subset \mathbb{R}^{2}$ as a bounded domain. For the ease of presentation, we consider mixed partial Neumann and partial Dirichlet boundary conditions in this paper. 
Specifically, the boundary $\partial\Omega$ is divided into the following:
$$
\partial\Omega=\Gamma_{\rm d}\cup\Gamma_{\rm t} \quad ~\mbox{and}~ \quad \partial\Omega=\Gamma_{\rm p}\cup\Gamma_{\rm f},
$$
where $\Gamma_{\rm d}$ and $\Gamma_{\rm t}$ are for displacement and stress boundary conditions; $\Gamma_{\rm p}$ and $\Gamma_{\rm f}$ are for pressure and flux boundary conditions. ~
Accordingly, the boundary conditions are the following: 
\begin{eqnarray} \label{BC_mix}
\ub &=& \bm{0} \quad \mbox{on} ~\Gamma_{\rm d}, ~~~~~~~(\sigma(\ub)-\alpha p \mathbf{I}) \cdot \nb = \mathbf{t} \quad \mbox{on} ~\Gamma_{\rm t}, \\
p &=& 0\quad \mbox{on} ~\Gamma_{\rm p}, ~~~~~~~\bm{z} \cdot \nb = g_2 \quad \mbox{on} ~\Gamma_{\rm f},
\end{eqnarray}
where $\sigma(\ub)$ is the deviatoric stress. For simplicity, the Dirichlet conditions are assumed to be homogeneous.

\section{Formulation of the Biot model as a saddle-point problem}


\subsection{Discrete formulation: $\mathbf{P}_1-\mathbf{P}_1-P_0$}
We apply the finite element method with domains normally shaped as triangles in $\mathbb{R}^{2}$. Let $\mathcal{T}_{h}$ be a partition of $\Omega$ into non-overlapping elements $K$. We denote by $h$ the size of the largest element in $\mathcal{T}_{h}$.
On the given partition $\mathcal{T}_{h}$  we apply the following finite element spaces \cite{Berger2015}:
%
%
\begin{eqnarray}
\bm{V}_h  &:=& \{ \bm{u}_h \in (C^0(\Omega))^2: \mathbf{u}_h|K \in {\bf P}_1(K) ~\forall K \in \mathcal{T}^h, \bm{u}_h = 0~ {\rm on}~ \Gamma_{\rm d} \},  \\
\bm{W}_h  &:=& \{ \bm{z}_h \in (C^0(\Omega))^2: \mathbf{z}_h|K \in {\bf P}_1(K) ~\forall K \in \mathcal{T}^h, \bm{z}_h \cdot \bm{n} = 0~ {\rm on}~ \Gamma_{\rm f} \},  \\
Q_h &:=& \{ p_h : p_h|K \in {\bf P}_0(K) ~ \forall K \in \mathcal{T}^h \} {\rm .}
\end{eqnarray}
The problem is to find $(\bm{u}_h^n, \bm{z}_h^n, p_h^n) \in \bm{V}_h \times \bm{W}_h \times Q_h$ ~~such that
\begin{equation}
\left\{
\begin{array}{l}
a(\bm{u}_h^n, \bm{v}_h)  -(p_h^n, \nabla \cdot \bm{v}_h) = (\bm{f}^n, \bm{v}_h) + (\bm{t}^n,  \bm{v}_h)_{\Gamma_{\rm t}}, ~\forall \bm{v}_h \in \bm{V}_h   \\
({\bf K}^{-1} \bm{z}_h^n, \bm{w}_h) - (p_h^n, \nabla \cdot \bm{w}_h) = (\bm{b}^n, \bm{w}_h), ~\forall \bm{w}_h \in \bm{W}_h      \\
(\nabla \cdot \bm{u}_{\Delta t,h}^n, q_h) + \frac{1}{\alpha}(\nabla \cdot \bm{z}_h^n, q_h)+\frac{c_0}{\alpha}(p_h^n, q_h)  + J(p_{\Delta t,h}^n,q_h) = \frac{1}{\alpha} (g^n, q_h), ~\forall q_h \in Q_h

\end{array}
\right.
\end{equation}
where
\begin{displaymath}
J(p,q) = \delta_{\rm STAB} \sum_K \int_{\partial K \backslash \partial \Omega} h_{\partial K} [p][q]ds
\end{displaymath}
is a stabilizing term \cite{Burman2007}, $p_{\Delta t,h}^n = (p_h^n - p^{n-1}_h)/\Delta t$, and $\bm{u}_{\Delta t,h}^n = (\bm{u}_h^n - \bm{u}^{n-1}_h)/\Delta t$.
%
%
The finite element discretization will lead to a twofold saddle-point problem of the following form:
\begin{equation}\label{matrix_A}
\left[\begin{array}{ccc}
A_{\bm{u}}  & 0      &  B^{T}_{1} \\[1mm]
0  & A_{\bm{z}}      & B^{T}_{2} \\[1mm]
B_1        &B_2    & -A_p \\
\end{array}\right]
\left[\begin{array}{c}
\bm{u}_h \\[1mm]
\bm{z}_h \\[1mm]
p_h
\end{array}\right]=
\left[\begin{array}{c}
{\bm f}_1 \\[1mm]
{\bm f}_2\\[1mm]
{\bm f}_3
\end{array}\right].
\end{equation}

\section{FETI-DP formulation for Biot model with discontinuous pressure field}

In the algorithm of FETI-DP \cite{Farhat2001}, the domain $\Omega$ is decomposed to $N$ nonoverlapping subdomains $\Omega_i$. Each subdomain is with the diameter in the order of $H$, and the neighboring subdomains are matched across the subdomain interface, $\Gamma = (\cup \partial \Omega_i) \setminus \partial \Omega$.

\subsection{FETI-DP algorithm for Biot model: interior and interface spaces}

We decompose the discrete displacement space $\bm{V}$, Darcy flow flux space $\bm{W}$ into interior and interface spaces $(\Vb = \Vb_{\rm I} \oplus \Vb_{\Gamma},~\Wb = \Wb_{\rm I} \oplus \Wb_{\Gamma})$.~
The discontinuous pressure space $Q$ is decomposed to the constant space $Q_{0}$ with the constant pressures on each subdomain and interior space $Q_{\rm I}$ which has the average zero over each subdomain.
Here $\bm{V}_{\rm I}$, $\bm{W}_{\rm I}$, and $Q_{\rm I}$ are the direct sums of subdomain interior spaces, and $\Vb_{\rm I} = \oplus_{i=1}^N \Vb_{\rm I}^i, ~\Wb_{\rm I} = \oplus_{i=1}^N \Wb_{\rm I}^i, ~Q_{\rm I} = \oplus_{i=1}^N Q_{\rm I}^i.$

\subsection{FETI-DP algorithm for Biot model: primal and dual variables}

The interface space $\Vb_{\Gamma}$ is further decomposed to primal and dual spaces:
\begin{equation}
\bm{V}_{\Gamma} = \bm{V}_{\Delta} \oplus \bm{V}_{\Pi} = (\oplus_{i=1}^N \bm{V}_{\Delta}^i) \oplus \bm{V}_{\Pi},
\end{equation}
where \noindent $\bm{V}_{\Pi}$ is the continuous, coarse level, and primal space. $\bm{V}_{\Delta}$ is the direct sum of independent subdomain dual interface spaces $\bm{V}_{\Delta}^i$ \cite{Tu2013}. Similarly $\bm{W}_{\Gamma}$ is decomposed to $\bm{W}_{\Delta}$ and $\bm{W}_{\Pi}$.\\\\
Let us represent $\bm{u}$ and $\bm{z}$ together as $\bm{U} = (\bm{u}, \bm{z}) \in \bm{V} \times \bm{W} $. The problem turns out to find $({\bm u}_{\rm I}, {\bm z}_{\rm I}, p_{\rm I}, {\bm u}_{\Pi}, {\bm z}_{\Pi}, {\bm u}_{\Delta}, {\bm z}_{\Delta}, p_0)$ $\in \bm{V}_{\rm I} \times \bm{W}_{\rm I} \times Q_{\rm I} \times \bm{V}_{\Pi} \times \bm{W}_{\Pi} \times \bm{V}_{\Delta} \times \bm{W}_{\Delta} \times Q_0$ such that
\begin{equation}
\left[\begin{array}{ccccc}
A_{\rm II} & B_{\rm II}^T &A_{\Pi \rm I}^T & A_{\Delta \rm I}^T & 0\\
B_{\rm II} & 0 & B_{{\rm I} \Pi} & B_{{\rm I} \Delta} & 0\\
A_{\Pi \rm I} & B_{{\rm I} \Pi}^T & A_{\Pi \Pi} & A_{\Delta \Pi}^T & B_{0 \Pi}^T\\
A_{\Delta \rm I} & B_{{\rm I} \Delta}^T & A_{\Delta \Pi} & A_{\Delta \Delta} & B_{0 \Delta}^T\\
0 & 0 & B_{0 \Pi} & B_{0 \Delta} & 0
\end{array}
\right]
\left[\begin{array}{c}
{\bm U}_{\rm I}    \\
p_{\rm I}    \\
{\bm U}_{\Pi}  \\
{\bm U}_{\Delta}\\
p_0  
\end{array}
\right]
=
\left[\begin{array}{c}
{\bm f}_{\rm I}    \\
0   \\
{\bm f}_{\Pi}  \\
{\bm f}_{\Delta} \\
0
\end{array}
\right].
\end{equation}

\subsection{FETI-DP algorithm for Biot model: Schur complement}

A Schur complement operator $\tilde{S}$ is defined in the following:
\begin{equation}
\left[\begin{array}{ccccc}
A_{\rm II} & B_{\rm II}^T &A_{\Pi \rm I}^T & 0 & A_{\Delta I}^T\\
B_{\rm II} & 0 & B_{{\rm I} \Pi} & 0 & B_{{\rm I} \Delta}\\
A_{\Pi \rm I} & B_{{\rm I} \Pi}^T & A_{\Pi \Pi} & B_{0 \Pi}^T & A_{\Delta \Pi}^T\\
0 & 0 & B_{0 \Pi} & 0 & B_{0 \Delta} \\
A_{\Delta \rm I} & B_{{\rm I} \Delta}^T & A_{\Delta \Pi} & B_{0 \Delta}^T & A_{\Delta \Delta}
\end{array}
\right]
\left[\begin{array}{c}
{\bm U}_{\rm I}    \\
p_{\rm I}    \\
{\bm U}_{\Pi}  \\
p_0\\
{\bm U}_{\Delta}
\end{array}
\right]
=
\left[\begin{array}{c}
0    \\
0   \\
0  \\
0 \\
\tilde{S}{\bm U}_{\Delta}  
\end{array}
\right].
\end{equation}

\noindent We introduce Lagrange multiplier $\lambda$ and the jump operator $B_{\Delta}$ to enforce the continuity of ${\bm U}_{\Delta}$ across $\Gamma$ \cite{Li2005}:
\begin{equation}
\left[\begin{array}{cc}
\tilde{S} & B_{\Delta}^T\\
B_{\Delta} & 0\\
\end{array}
\right]
\left[\begin{array}{c}
{\bm U}_{\Delta} \\
{\bm \lambda}
\end{array}
\right]
=
\left[\begin{array}{c}
{\bm f}_{\Delta}^*    \\
0
\end{array}
\right].
\end{equation}
The problem is reduced to find $\lambda \in \Lambda = B_\Delta \bm{U}_{\Delta}$ such that
\begin{equation}
B_{\Delta} \tilde{S}^{-1}B_{\Delta}^T \lambda = B_{\Delta} \tilde{S}^{-1} {\bm f}_{\Delta}^*.
\end{equation}

\noindent This is solved by Preconditioned Conjugate Gradient (PCG).

\subsection{Dirichlet preconditioner}

We define a Schur complement operator, the discrete Harmonic $H_{\Delta}^{(i)}$ on $\Omega_i$ as follows:
\begin{equation} \label{harmonic_extension}
\left[\begin{array}{cc}
A_{\rm II}^{(i)} & A_{{\rm I} \Delta}^{(i)} \\
A_{\Delta I}^{(i)} & A_{\Delta \Delta}^{(i)}
\end{array}
\right]
\left[\begin{array}{c}
{\bm U}_{\rm I}^{(i)} \\
{\bm U}_{\Delta}^{(i)}
\end{array}
\right]
=
\left[\begin{array}{c}
0   \\
H_{\Delta}^{(i)} {\bm U}_{\Delta}^{(i)}
\end{array}
\right].
\end{equation}

\noindent The Dirichlet preconditioner is formulated in the following:
\begin{equation}
M_{\lambda, D}^{-1} = B_{\Delta, D} H_{\Delta} B_{\Delta, D}^T,
\end{equation}
where $B_{\Delta, D}$ is a scaled operator obtained from $B_{\Delta}$ by the scaling factor $1/N_x$ with $N_x$ as the number of subdomains sharing each node $x$ in the interface $\Gamma$. $H_{\Delta}$ is the direct sum of $H_{\Delta}^{(i)}$ \cite{Tu2015}.

\section{Numerical experiments}
%
A test problem is formulated with $\alpha = 1$, $c_0 = 0$, $\Omega = [0,1]^2$, and $t \in [0, 0.25]$:
\begin{eqnarray}
-(\lambda + \mu) \nabla (\nabla \cdot \bm{u}) - \mu \nabla^2 \bm{u} + \nabla p &=& 0,  \label{test1} \nonumber \\
 \mathbf{K}^{-1} \bm{z} + \nabla p &=& 0,   \label{test2}\\
\nabla \cdot (\bm{u}_t + \bm{z}) &=& g_1. \label{test3} \nonumber
\end{eqnarray}
The involving initial and boundary conditions are the following:
\begin{equation}
\left\{
\begin{array}{l}
\ub~=\bm{0}\quad \mbox{on} ~\partial \Omega = \Gamma_{\rm d}, \\
\bm{z} \cdot\nb~=g_2 \quad \mbox{on} ~\partial \Omega = \Gamma_{\rm f},\\
\bm{u}(\bm{x},0) = 0, \bm{x} \in \Omega,\\
p(\bm{x},0) = 0, \bm{x} \in \Omega.
\end{array}
\right.
\end{equation}
%
%
We consider the following analytic solution:
\begin{eqnarray}
\bm{u}  &=& \frac{-1}{4 \pi (\lambda + 2 \mu)} \left[\begin{array}{c}
\cos(2 \pi x) \sin(2 \pi y) \sin(2 \pi t)    \\
\sin(2 \pi x) \cos(2 \pi y) \sin(2 \pi t)
\end{array}
\right],  \nonumber \\
\bm{z}  &=& -2 \pi k \left[\begin{array}{c}
\cos(2 \pi x) \sin(2 \pi y) \sin(2 \pi t)    \\
\sin(2 \pi x) \cos(2 \pi y) \sin(2 \pi t)
\end{array}
\right], \\
p &=& \sin(2 \pi x) \sin(2 \pi y) \sin(2 \pi t), \nonumber
\end{eqnarray}
and derive the compatible source term of $g_1$.
%

\subsection{Numerical implementation}
In the implementation of finite elements, we use a finite element library, libMesh \cite{Kirk2006}. We apply triangular elements with 3 nodes. For domain partitioning, we apply ParMETIS \cite{Karypis2011}. Krylov subspace iterative main solver of Preconditioned Conjugate Gradient (PCG) and FETI-DP algorithms are based on PETSc \cite{Balay1997, Balay2021_manual, Balay2021_PETSc_webpage} and KSPFETIDP and PCBDDC classes within PETSc \cite{Zampini2016}. The initial guess is zero and the stopping criterion is set to be $10^{-8}$, the reduction of the residual norm. The stabilizing factor is $\delta_{\rm STAB}=100$, the time-stepping is $dt = 0.00625$, and Young's modulus $E = 1000$ Pa. In each test, we count the iteration of the FETI-DP solver.
\subsection{Scalability of FETI-DP algorithms} 
Scalability of FETI-DP preconditioning for the Biot model is tested with increasing number of subdomains $N$. The subdomain size $H/h$ is set with 8, 12, or 16. In the first case, $(\nu = 0.3, k = 10^{-2})$ of compressible elasticity and permeable Darcy flow is tested with Dirichlet preconditioner. As shown in Table 1, FETI-DP iteration numbers are bounded when subdomains were increased from $2 \times 2$ to $8 \times 8$. In the second case, $(\nu = 0.4999, k = 10^{-7})$ of almost incompressible elasticity and less permeable Darcy flow is tested with Dirichlet preconditioner. FETI-DP iteration numbers are larger than the first case, but still bounded while subdomains are increased from $2 \times 2$ to $8 \times 8$, showing no issues of elastic locking. This is consistent with a theoretical scalability of FETI-DP for almost incompressible elasticity \cite{Klawonn2005}. 

\begin{table}[bthp]
\begin{center}
  \centering
  \small
  \begin{tabular}{c|cc|cc|cc}
    \hline
&\multicolumn{2}{c|}{ $H/h=8$} & \multicolumn{2}{c|}{ $H/h=12$} & \multicolumn{2}{c}{ $H/h=16$}\\
    \hline      
&\multicolumn{1}{c|}{ $\nu = 0.3$} & \multicolumn{1}{c|}{ $\nu = 0.4999$} &\multicolumn{1}{c|}{ $\nu = 0.3$} & \multicolumn{1}{c|}{ $\nu = 0.4999$} &\multicolumn{1}{c|}{ $\nu = 0.3$} & \multicolumn{1}{c}{ $\nu = 0.4999$}\\    
&\multicolumn{1}{c|}{ $k = 10^{-2}$} & \multicolumn{1}{c|}{ $k = 10^{-7}$} &\multicolumn{1}{c|}{ $k = 10^{-2}$} & \multicolumn{1}{c|}{ $k = 10^{-7}$} &\multicolumn{1}{c|}{ $k = 10^{-2}$} & \multicolumn{1}{c}{ $k = 10^{-7}$}\\
    \hline  
         {$N$}         & \multicolumn{1}{c|}{iteration} & \multicolumn{1}{c|}{iteration} & \multicolumn{1}{c|}{iteration} & \multicolumn{1}{c|}{iteration} & \multicolumn{1}{c|}{iteration} & \multicolumn{1}{c}{iteration} \\

    \hline
    $2\times2$ &  4  &  9  &  4 &  10  & 4  & 11\\
    $3\times3$ &  5  &  9  &  5 &  12  & 5  & 14\\
    $4\times4$ &  5  &  9  &  5 &  12  & 5  & 17\\
    $5\times5$ &  5  &  10  &  5 & 13 & 5  & 15 \\
    $6\times6$ &  5  &  11  &  5 & 13 & 5  & 16\\
    $7\times7$ &  5  &  12  &  5 & 15 & 5  & 16\\
    $8\times8$ &  5  &  12  &  5 & 13 & 5  & 17\\
    \hline
  \end{tabular}
\caption{Scalability of the FETI-DP algorithms with Dirichlet preconditioner for the saddle-point problem of Biot model. Iteration counts for increasing number of subdomains $N$. Fixed $\delta_{\rm STAB} = 100, dt = 0.00625$, and $E=1000$.
}
\end{center}
\end{table}

\section{Conclusion}
We have explored the scalability of the FETI-DP algorithms for the 2D Biot model. Upon numerical scalabilities of compressible elasticity with Darcy's flow as well as almost incompressible elasticity with limited Darcy's flow, it remains to test parameter robustness, possibly in the presence of heterogeneity of parameters. Overall, the numerical results are a foundation for further advancement of scalable FETI-DP / BDDC preconditioners for poroelastic large deformation.
\begin{acknowledgement}
The author gives thanks to Dr. Mingchao Cai for the introduction of Biot models and invaluable discussions. The author was partly supported by NSF DMS-1831950, and the virtual attendance to DD27 conference by Penn State University NSF travel fund.
\end{acknowledgement}

\end{document}